\documentclass{svjour3}
\smartqed  
\usepackage{graphicx}
\usepackage{mathptmx}      
\usepackage{amsmath,amsfonts,amssymb}
\usepackage{bbold}
\newcommand\1{\mathbb{1}}
\newcommand\RR{\mathbb{R}}
\newcommand\NN{\mathbb{N}}
\newcommand\KK{\mathbb{K}}

\newcommand\Pb{\mathbf{P}}
\newcommand\Ex{\mathbf{E}}
\newcommand\dd{\mathrm{d}}
\newcommand\Bb{\mathcal{B}}
\newcommand\Cc{\mathcal{C}}
\newcommand\Dd{\mathcal{D}}
\newcommand\Ff{\mathcal{F}}
\newcommand\Hh{\mathcal{H}}
\newcommand\Kk{\mathcal{K}}
\newcommand\Ll{\mathcal{L}}
\newcommand\Nn{\mathcal{N}}
\newcommand{\eqd}{\mathrel{\overset{\textrm{d}}{=}}}
\newcommand{\eqdef}{\mathrel{\overset{\mathrm{\Delta}}{=}}}
\DeclareMathOperator*{\Liminf}{\underline{\lim}}
\newenvironment{lemproof}[2][]{\let\proofnameorig\proofname
  \renewcommand\proofname{Proof of Lemma~\ref{#2}#1}\begin{proof}}{\qed
  \end{proof}\let\proofname\proofnameorig}
\def\zzz#1.{\phantom{A}}
\journalname{\zzz}

\begin{document}

\title{On a Poissonian Change-Point Model with Variable Jump Size}

\subtitle{}

\titlerunning{}

\author{Sergue\"{\i} Dachian
        \and
        Lin Yang}

\authorrunning{}

\institute{S. Dachian \at
              Laboratoire de Math\'ematiques UMR6620, Universit\'e Blaise
              Pascal, Clermont-Ferrand, France\\
              \email{Serguei.Dachian@math.univ-bpclermont.fr}
           \and
           L. Yang \at
              Laboratoire Manceau de Math\'ematiques, Universit\'e du Maine,
              Le Mans, France\\
              \email{Lin\_Yang.Etu@univ-lemans.fr}}

\date{Received: / Accepted:}

\maketitle

\begin{abstract}
A model of Poissonian observation having a jump (change-point) in the
intensity function is considered.  Two cases are studied.  The first one
corresponds to the situation when the jump size converges to a non-zero limit,
while in the second one the limit is zero.  The limiting likelihood ratios in
these two cases are quite different.  In the first case, like in the case of a
fixed jump size, the normalized likelihood ratio converges to a log Poisson
process.  In the second case, the normalized likelihood ratio converges to a
log Wiener process, and so, the statistical problems of parameter estimation
and hypotheses testing are asymptotically equivalent in this case to the well
known problems of change-point estimation and testing for the model of a
signal in white Gaussian noise.  The properties of the maximum likelihood and
Bayesian estimators, as well as those of the general likelihood ratio, Wald's
and Bayesian tests are deduced form the convergence of normalized likelihood
ratios.  The convergence of the moments of the estimators is also established.
The obtained theoretical results are illustrated by numerical simulations.
\end{abstract}

\section{Introduction}

In regular statistical experiments, the limit of the normalized likelihood
ratio is always the same, because the families are LAN $\bigl($see, for
example,~\cite{IH81}$\bigr)$. In the case of non-regular statistical models
for Poisson processes, there exists a large diversity of limiting likelihood
ratio processes: change-point type models lead to a log Poisson process,
``cusp'' type singularities provide a log fBm process, while in the models
with $0$-type or $\infty$-type singularities the limit processes are more
sophisticated $\bigl($see, respectively,~\cite{Kut98,Da03,Da11}$\bigr)$. Note
that in change-point type models for diffusion processes, and particularly in
the model of a discontinuous signal in white Gaussian noise (WGN), the
limiting likelihood ratio is a log Wiener process $\bigl($see, for
example,~\cite{IH81,Kut04}$\bigr)$. It is interesting to investigate the
relations between the different limit processes. This study was initiated in
the recent works~\cite{Da10,DN}. The present work is a part of this
investigation, since we study a change-point model with variable jump size for
Poissonian observations, and we obtain two different limits depending on the
way the jump size is varying.

More precisely, we consider two cases. The first one corresponds to the
situation when the jump size converges to a non-zero limit, while in the
second one the limit is zero. The limiting likelihood ratios in these two
cases are quite different. In the first case, as one could expect, the
normalized likelihood ratio converges to a log Poisson process, just like the
case of a fixed jump size. In the second case, the normalized likelihood ratio
converges to a log Wiener process, that is, the statistical problems of
parameter estimation and hypotheses testing are asymptotically equivalent to
the well known problems of change-point estimation and testing for signal in
WGN model. Let us note, that even if the latter result may seem unexpected, it
is quite natural in the light of the recent work~\cite{Da10} of one of the
authors, where a relation between the log Poisson and the log Wiener limiting
likelihood ratios was discovered.

Let us also mention that this situation is somewhat similar to what happens in
the case of multi-phase regression models, where the limiting likelihood ratio
is a log compound Poisson process in the case of a fixed jump size, while it
is a log Wiener process in the case of a variable jump size converging to zero
$\bigl($see, for example,~\cite{Fuj} and the references therein$\bigr)$. Note
also, that the recent work~\cite{DN} shades the light on the latter case, just
as~\cite{Da10} do in our case of Poissonian observations.

Note finally, that we show not only the convergence of normalized likelihood
ratios, but also the convergence of the moments of the estimators. This last
convergence allows one, for example, to approximate the limiting mean square
errors of the maximum likelihood and Bayesian estimators in the case of
Poisson observations by the well known limiting mean square errors of these
estimators calculated for signal in WGN model.

The paper is organized as follows.  In Section~\ref{Sec1} we describe the
model of observations.  In Section~\ref{Sec2} we study the asymptotic behavior
of the likelihood ratio.  In Section~\ref{Sec3}, using the convergence of
normalized likelihood ratio obtained in Section~\ref{Sec2}, we study the
problem of parameter estimation.  Similarly, in Section~\ref{Sec4} we study
the problem of hypothesis testing and illustrate the results by numerical
simulations.  Finally, Section~\ref{Sec5} contains the proofs of all the
lemmas.

\section{Change-point model with variable jump size}
\label{Sec1}

Suppose we observe $n$ independent realizations
$X^{(n)}_j=\bigl\{X^{(n)}_j(t),\ t\in\left[0,\tau\right]\bigr\}$,
$j=1,\ldots,n$, of an inhomogeneous Poisson process on the interval
$\left[0,\tau\right]$ (the constant $\tau>0$ is supposed to be known) of
intensity measure
$$
\Lambda^{(n)}_\vartheta\bigl(A\bigr)=\int_A \lambda^{(n)}_\vartheta(t)\,\dd t,
\qquad A\in\Bb\bigl(\left[0,\tau\right]\bigr),
$$
with intensity function $\lambda^{(n)}_\vartheta$, where
$\vartheta\in\Theta=\left(\alpha,\beta\right)$, $0\leq\alpha<\beta\leq\tau$,
is some unknown parameter. The observation will be denoted
$X^{(n)}=\bigl\{X^{(n)}_1,\ldots,X^{(n)}_n\bigr\}$ and the corresponding
probability distribution will be denoted $\Pb^{(n)}_\vartheta$.

Let us note that this model of observation is equivalent to observing a single
realization on the interval $\left[0,n\tau\right]$ of an inhomogeneous Poisson
process with the $\tau$-periodic intensity function coinciding with
$\lambda^{(n)}_\vartheta$ on $\left[0,\tau\right]$.

The parameter $\vartheta$ corresponds to the location of a jump in the
(elsewhere continuous) intensity function~$\lambda^{(n)}_\vartheta$. The size
of the jump (depending on $n$) will be denoted $r_n$ and will be supposed
converging to some $r\in\RR$.  As we will see below, the behavior of our model
depends on either one has $r\ne 0$ or $r=0$ and is quite different in these
two cases.

More precisely, we assume that the following conditions are satisfied.

\medskip

\noindent\textbf{(C1)} The intensity function $\lambda^{(n)}_\vartheta(t)$ can be
written as $\lambda^{(n)}_\vartheta(t)=\psi_n(t)+r_n\1_{\{t>\vartheta\}}$, where the
function $\psi_n$ is continuous on $\left[0,\tau\right]$.

\medskip

\noindent\textbf{(C2)} For all $t\in\left[0,\tau\right]$, there exist the
$\lim\limits_{n \to +\infty}\psi_n(t)=\psi(t)>0$ and, moreover, this
convergence is uniform with respect to~$t$.

\medskip

\noindent\textbf{(C3)} As $n\to+\infty$, the jump size $r_n$ converges
to some $r\in\RR$, that is, $r_n\to r$.  In the case $r=0$, we also
suppose that this convergence ($r_n\to0$) is slower than $n^{-1/2}$,
that is, $n\,r_n^2\to+\infty$.

\medskip

\noindent\textbf{(C4)} The family of functions
$\bigl\{\lambda_\vartheta^{(n)}\bigr\}_{n\in\NN,\vartheta\in\Theta}$ is uniformly
strictly positive and uniformly bounded, that is, there exist some constants
$\ell,L>0$ such that
$$
\ell\leq\lambda_\vartheta^{(n)}(t)\leq L
$$
for all $n\in\NN$, $\vartheta\in\Theta$ and $t\in[0,\tau]$.

\medskip

Note that the conditions \textbf{C1}~--~\textbf{C3}, together with the natural
condition
\begin{equation}
\label{I4bis}
r>-\min_{t\in\left[0,\tau\right]}\psi(t),
\end{equation}
easily imply that the condition~\textbf{C4} holds for the family
$\bigl\{\lambda_\vartheta^{(n)}\bigr\}_{n\geq n_0,\vartheta\in\Theta}$ with
some $n_0\in\NN$. So, in the asymptotic setting ($n\to+\infty$), the
condition~\textbf{C4} can be replaced by~\eqref{I4bis}, and we
assume~\textbf{C4} instead of the latter only for convenience (as well as in
order for our model to be well defined for all $n\in\NN$). Note also that in
the case $r=0$, the condition~\eqref{I4bis} is automatically satisfied.

An important particular case of this model is when only the jump size
$\bigl($and not the regular part of $\lambda^{(n)}_\vartheta\bigr)$ depend on
$n$. More precisely, the conditions \textbf{C1}~--~\textbf{C2} will be clearly
met if we assume that the following condition is satisfied.

\bigskip

\noindent\textbf{(C0)} The intensity function $\lambda^{(n)}_\vartheta(t)$ can
be written as $\lambda^{(n)}_\vartheta(t)=\psi(t)+r_n\1_{\{t>\vartheta\}}$,
where the function $\psi$ is strictly positive and continuous on
$\left[0,\tau\right]$.

\section{Asymptotic behavior of the likelihood ratio}
\label{Sec2}

The likelihood of our model is given by $\bigl($see, for
example,~\cite{Kut84}$\bigr)$
\begin{equation}
\label{LR}
\begin{aligned}
L_n\bigl(\vartheta,X^{(n)}\bigr)
&=\exp\biggl\{\sum_{j=1}^n\int_{\left[0,\tau\right]}\ln
\lambda_\vartheta^{(n)}(t)\;X^{(n)}_j(\dd t) - n\int_0^\tau
\bigl[\lambda_\vartheta^{(n)}(t)-1\bigr]\,\dd t\biggr\}\\
&=\exp\biggl\{\sum_{j=1}^n \ \sum_{i\in
  I_j^{(n)}}\ln\lambda_\vartheta^{(n)}(t_{j,i}) - n\int_0^\tau
\bigl[\lambda_\vartheta^{(n)}(t)-1\bigr]\,\dd t\biggr\},
\end{aligned}
\end{equation}
where $t_{j,i}$, $i\in I_j^{(n)}$, are the jump times of the
process~$X^{(n)}_j$.  Note that as function of $\vartheta$, each
$\lambda_\vartheta^{(n)}(t_{j,i})$ is discontinuous (has a jump and is right
continuous) at $\vartheta=t_{j,i}$.  So, $L_n\bigl({\,\cdot\,},X^{(n)}\bigr)$
is a random process with c\`adl\`ag (continuous from the right and having
finite limits from the left) trajectories.

We put $\varphi_n=\frac1n$ in the case $r\ne0$ and
$\varphi_n=\frac{1}{n\,r_n^2}$ in the case $r=0$, and we introduce the
normalized likelihood ratio
\begin{align*}
Z_{n,\vartheta}(u)&=\frac{L_n\bigl(\vartheta+u\varphi_n,
  X^{(n)}\bigr)}{L_n\bigl(\vartheta,X^{(n)}\bigr)}\\
&=\exp\biggl\{\sum_{j=1}^n\int_{\left[0,\tau\right]}\ln
\frac{\lambda_{\vartheta+u\varphi_n}^{(n)}(
  t)}{\lambda_\vartheta^{(n)}(t)}\,X^{(n)}_j(\dd t)
-n\int_0^{\tau}\bigl(\lambda_{\vartheta+u\varphi_n}^{(n)}(t)
-\lambda_\vartheta^{(n)}(t)\bigr)\,\dd t\biggr\}\\
&=\exp\biggl\{\sum_{j=1}^n\ \sum_{i\in
  I_j^{(n)}}\ln\frac{\lambda_{\vartheta+u\varphi_n}^{(n)}(
  t_{j,i})}{\lambda_\vartheta^{(n)}(t_{j,i})}
-n\int_0^{\tau}\bigl(\lambda_{\vartheta+u\varphi_n}^{(n)}(t)
-\lambda_\vartheta^{(n)}(t)\bigr)\,\dd t\biggr\},
\end{align*}
where $u\in U_n=\bigl(\varphi_n^{-1}(\alpha-\vartheta),\varphi_n^{-1}
(\beta-\vartheta)\bigr)$.

Note that in both cases we have (by the condition~\textbf{C3} in the case
$r=0$) $\varphi_n\to 0$.

Note also that if $u>0$, we can rewrite $Z_{n,\vartheta}(u)$ as
\begin{equation}
\label{explicitSI1}
\begin{aligned}
Z_{n,\vartheta}(u)&=\exp\biggl\{\sum_{j=1}^n\int_{\left(\vartheta,\vartheta+u\varphi_n\right]}\ln
  \frac{\psi_n(t)}{\psi_n(t)+r_n}\;X^{(n)}_j(\dd t)
+n\int_\vartheta^{\vartheta+u\varphi_n}r_n\,\dd t\biggr\}\\
&=\exp\biggl\{\sum_{j=1}^n\ \sum_i\ln
\frac{\psi_n(t_{j,i})}{\psi_n(t_{j,i})+r_n}+u\,r_n^\gamma\biggr\}.
\end{aligned}
\end{equation}
In the last expression the inner sum is taken over the set $\bigl\{i\in
I_j^{(n)}\::\:\vartheta<t_{j,i}\leq\vartheta+u\varphi_n\bigr\}$ and we have
$\gamma=1$ in the case $r\ne 0$ and $\gamma=-1$ in the case $r=0$.

Similarly, if $u<0$, we have
\begin{align*}
Z_{n,\vartheta}(u)&=\exp\biggl\{\sum_{j=1}^n\int_{\left(\vartheta+u\varphi_n,\vartheta\right]}\ln
  \frac{\psi_n(t)+r_n}{\psi_n(t)}\;X^{(n)}_j(\dd t)
-n\int_{\vartheta+u\varphi_n}^\vartheta r_n\,\dd t\biggr\}\\
&=\exp\biggl\{\sum_{j=1}^n\ \sum_i\ln
\frac{\psi_n(t_{j,i})+r_n}{\psi_n(t_{j,i})}+u\,r_n^\gamma\biggr\},
\end{align*}
where the last sum is taken over the set $\bigl\{i\in
I_j^{(n)}\::\:\vartheta+u\varphi_n<t_{j,i}\leq\vartheta\bigr\}$ and $\gamma$
is as above.

Note equally, that the process $\ln Z_{n,\vartheta}$ has independent
increments. Indeed, its increments on disjoint intervals involve stochastic
integrals (of a deterministic function with respect to Poisson processes) on
disjoint intervals, and hence are independent. In other words, using the
terminology of Strasser~\cite{Str85}, our model is an ``experiment with
independent increments''. Note also, that in this case the process
$Z_{n,\vartheta}$ $\bigr($as well as, for example, the process
$Z_{n,\vartheta}^{1/2}\bigr)$ is clearly a Markov process.

Note finally, that the trajectories of the process $Z_{n,\vartheta}$ are
c\`adl\`ag functions.  Moreover, correctly extending these trajectories to the
whole real line, one can consider that they belong to the Skorohod space
$\Dd_0(\RR)$.  This space is defined as the space of functions $f$ on $\RR$
which do not have discontinuities of the second kind and which are vanishing
at infinity, that is, such that
$\lim\limits_{u\to\pm\infty}f(u)=0$. We assume that all the functions
$f\in\Dd_0(\RR)$ are continuous from the right (are c\`adl\`ag).

Let us recall that the Skorohod metric on the space $\Dd_0(\RR)$ is introduced
by
$$
d(f,g)=\inf_\lambda\biggl[\sup_{u\in\RR}
  \left|f(u)-g\bigl(\lambda(u)\bigr)\right|+
  \sup_{u\in\RR}\left|u-\lambda(u)\right|\biggr],
$$
where the $\inf$ is taken over all strictly increasing continuous one-to-one
mappings $\lambda:\RR \to \RR$.

Let us also recall a criterion of weak convergence in $\Dd_0(\RR)$.  We put
$$
\Delta_h(f)=\sup_{u\in\RR}\
\sup_{u',u''}\Bigl[\min\bigl\{\left|f(u') -
  f(u)\right|,\left|f(u'')-f(u)\right|\bigr\}\Bigr]+\sup_{|u|>1/h}|f(u)|,
$$
where the inner $\sup$ is over all $u',u''$ such that $u-h \leq u'<u \leq
u''<u+h$.  A criterion of weak convergence in $\Dd_0(\RR)$ is given in the
following lemma $\bigl($see~\cite{Skoro69} for more details$\bigr)$.

\begin{lemma}
\label{CWC}
Let\/ $z_{n,\vartheta}$, $n\in\NN$, and\/ $z_\vartheta$ be random processes with
realizations belonging to\/ $\Dd_0(\RR)$ with probability\/ $1$. If, as\/
$n\to+\infty$, the finite dimensional distributions of\/
$z_{n,\vartheta}$ converge uniformly in\/ $\vartheta\in\KK$ to the finite
dimensional distributions of\/ $z_\vartheta$, and if for any\/ $\delta>0$
\begin{equation}
\label{Skoro}
\lim_{h\to0}\ \sup\limits_{n\in\NN,\vartheta\in\KK}\,\Pb
\bigl\{\Delta_h(z_{n,\vartheta})>\delta\bigr\}=0,
\end{equation}
then, uniformly in\/ $\vartheta\in\KK$, the process\/ $z_{n,\vartheta}$ converges
weakly in the space\/ $\Dd_0(\RR)$ to the process\/~$z_\vartheta$.
\end{lemma}

Note that here and in the sequel $\KK$ denotes an arbitrary compact in
$\Theta$.

The main objective of this section is the study of the asymptotic behavior
$\bigl($in the sense of the weak convergence in the space $\Dd_0(\RR)$ as
$n\to\infty\bigr)$ of the above introduced normalized likelihood ratio
$Z_{n,\vartheta}$.  This behavior depends on either one has $r\ne 0$ or $r=0$
and is quite different in these two cases, so the limit process must be
introduced in a different manner in these two cases.

\paragraph{Case $r\ne 0$ limit process}
In the case $r\ne 0$, the limit process is a log Poisson type process and is
introduced by
$$
Z_\vartheta(u)=\begin{cases}
\vphantom{\displaystyle\int}
\exp\Bigl\{\ln\frac{\psi(\vartheta)}{\psi(\vartheta)+r}\,X^+(u) + r
u\Bigr\},&\text{if }u\geq0,\\
\vphantom{\displaystyle\int}
\exp\Bigl\{\ln\frac{\psi(\vartheta)+r}{\psi(\vartheta)}\,X^-\bigl((-u)-\bigr)
+ r u\Bigr\},&\text{if }u<0,
\end{cases}
$$
where $X^+$ and $X^-$ are independent Poisson processes on $\RR_+$
of constant intensities $\psi(\vartheta)+r$ and $\psi(\vartheta)$
respectively.

Let us note that $Z_\vartheta(u)\eqd Z^*_\rho(-r u)$ with the constant
$\rho=\bigl|\ln\frac{\psi(\vartheta)}{\psi(\vartheta)+r}\bigr|$ and the
process $Z^*_\rho$ defined by
$$
Z_\rho^*(v)=\begin{cases}
\exp\bigl\{\rho Y^+(v)-v\bigr\},&\text{if }v\geq0,\\
\exp\bigl\{-\rho Y^-\bigl((-v)-\bigr)-v\bigr\},&\text{if }v<0,
\end{cases}
$$
where $Y^+$ and $Y^-$ are independent Poisson processes on $\RR_+$ of constant
intensities $\frac1{e^\rho-1}$ and $\frac1{1-e^{-\rho}}$ respectively.

Note also that the process $Z^*_\rho$ was recently studied in~\cite{Da10} and
that its trajectories (as well as those of the process~$Z_\vartheta$) almost
surely belong to the space $\Dd_0(\RR)$.  $\bigl($More rigorously, in order to
keep all the trajectories in the space $\Dd_0(\RR)$, above we should rather
have written $Z_\vartheta(u)\eqd Z^*_\rho\bigl((-r u)-\bigr)$ in the case
$r>0\bigr)$.

\paragraph{Case $r=0$ limit process}
In the case $r=0$, the limit process is a log Wiener type process and is
introduced by
$$
Z_\vartheta(u)=
\exp\biggl\{\psi^{-1/2}(\vartheta)\,W(u)-\frac{|u|}{2\psi(\vartheta)}\biggr\},\quad
u\in\RR,
$$
where $W(u)$, $u\in\RR$, is a double-sided Brownian motion (Wiener process).

Let us note that $Z_\vartheta(u)\eqd Z^*\bigl(u/\psi(\vartheta)\bigr)$ with
the process $Z^*$ defined by
\begin{equation}
\label{Zstar}
Z^*(v)=\exp\biggl\{W(v)-\frac{|v|}{2}\biggr\},\quad v\in\RR.
\end{equation}

Note also that the trajectories of the processes $Z^*$ and $Z_\vartheta$
almost surely belong to the space $\Cc_0(\RR)$ of continuous functions on
$\RR$ vanishing at infinity, and that $\Cc_0(\RR)\subset\Dd_0(\RR)$.

\paragraph{}

Now we can state the following theorem about the asymptotic behavior of the
normalized likelihood ratio.

\begin{theorem}
\label{CLR}
Let the conditions\/ \emph{\textbf{C1}~--~\textbf{C4}} be fulfilled. Then,
uniformly in\/ $\vartheta\in\KK$, the process\/ $Z_{n,\vartheta}$ converges
weakly in the space\/ $\Dd_0(\RR)$ to the process\/ $Z_\vartheta$.
\end{theorem}

Let us also remark, that sometimes it may be more convenient to use a slightly
different rate for introducing the normalized likelihood ratio.  More
precisely, one can use the rate $\varphi_n^*=\frac1{\left|r\right|n}$ (rather
than $\varphi_n=\frac1n$) in the case $r\ne 0$, and the rate
$\varphi_n^*=\frac{\psi(\vartheta)}{n\,r_n^2}$ (rather than
$\varphi_n=\frac1{n\,r_n^2}$) in the case $r=0$.  That is, one can consider
(instead of $Z_{n,\vartheta}$) the normalized likelihood ratio
$Z^*_{n,\vartheta}$ defined by
$$
Z^*_{n,\vartheta}(v)=\frac{L_n\bigl(\vartheta+v\,\varphi_n^*,
  X^{(n)}\bigr)}{L_n\bigl(\vartheta,X^{(n)}\bigr)}=
Z_{n,\vartheta}\bigl(c\,v\bigr)
$$
with $c=1/\left|r\right|$ in the case $r\ne 0$, and $c=\psi(\theta)$ in the
case $r=0$.  Then, Theorem~\ref{CLR} will be clearly transformed to the
following (equivalent) statement.

\begin{theorem}
Let the conditions\/ \emph{\textbf{C1}~--~\textbf{C4}} be fulfilled. Then,
uniformly in\/ $\vartheta\in\KK$, the process\/ $Z^*_{n,\vartheta}$ converges
weakly in the space\/ $\Dd_0(\RR)$ to
\begin{itemize}
\item the process\/ $Z^*_\rho$, in the case\/ $r<0$,
\item the process\/ $Z^\star_\rho$ defined by\/
  $Z^\star_\rho(v)=Z^*_\rho\bigl((-v)-\bigr)$, in the case\/ $r>0$,
\item the process\/ $Z^*$, in the case\/ $r=0$.
\end{itemize}
\end{theorem}

The proof of Theorem~\ref{CLR} consist in checking the criterion of week
convergence given in Lemma~\ref{CWC}.  For this, we follow the methods and
ideas used in~\cite[Chapters~5.3 and~5.4]{IH81} and establish several lemmas
(the proofs of the lemmas are in Section~\ref{Sec5}).

\begin{lemma}
\label{L2}
Let the conditions\/ \emph{\textbf{C1}~--~\textbf{C4}} be fulfilled. Then the
finite-dimensional distributions of the process\/ $Z_{n,\vartheta}$ converge to
those of the process\/ $Z_\vartheta$, and this convergence is uniform with
respect to\/ $\vartheta\in\KK$.
\end{lemma}

\begin{lemma}
\label{L3}
Let the conditions\/ \emph{\textbf{C1}~--~\textbf{C4}} be fulfilled. Then
there exists a constant\/ $C>0$ such that
$$
\Ex^{(n)}_{\vartheta}\bigl|Z_{n,\vartheta}^{1/2}(u_1)-Z_{n,\vartheta}^{1/2}(u_2)\bigr|^2\leq
C\left|u_1-u_2\right|
$$
for all\/ $n\in\NN$, $u_1,u_2\in U_n$ and\/ $\vartheta\in\KK$.
\end{lemma}

\begin{lemma}
\label{L4}
Let the conditions\/ \emph{\textbf{C1}~--~\textbf{C4}} be fulfilled. Then
there exists a constant\/ $k_*>0$ such that
$$
\Ex^{(n)}_{\vartheta}Z_{n,\vartheta}^{1/2}(u)\leq\exp\bigl\{-k_*\left|u\right|\bigr\}
$$
for all\/ $u\in U_n$, $\vartheta\in\KK$ and sufficiently large values of\/ $n$
(all\/ $n\in\NN$ in the case $r=0$).
\end{lemma}

\paragraph{Final argument of the proof of Theorem~\ref{CLR} in the case $r\ne
  0$}
In this case, defining $Z_{n,\vartheta;\mathrm{a.c.}}^{1/2}$ to be the
absolutely continuous component of the function $Z_{n,\vartheta}^{1/2}$ and,
for $p=1,2$, denoting $A_p=A_p(u,u+h)$ the event that $Z_{n,\vartheta}$ has at
least $p$ jumps on the interval $\left(u,u+h\right)$, we also have the
following lemma.

\begin{lemma}
\label{L6}
Let the conditions\/ \emph{\textbf{C1}~--~\textbf{C4}} be fulfilled with\/
$r\ne 0$. Then the inequalities
\begin{align}
\Ex^{(n)}_{\vartheta}\bigl|Z_{n,\vartheta;\mathrm{a.c.}}^{1/2}(u+h)-
Z_{n,\vartheta;\mathrm{a.c.}}^{1/2}(u)\bigr|^2&\leq C h^2,\notag\\
\label{dis_1}
\Pb_\vartheta^{(n)}(A_1)&\leq D_1 h\\
\intertext{and}
\label{dis_2}
\Pb_\vartheta^{(n)}(A_2)&\leq D_2h^2
\end{align}
hold with certain constants\/ $C,D_1,D_2>0$ (independent of\/ $n$,
$\vartheta$, $u$ and\/ $h$).
\end{lemma}

Now, with the help of the above lemmas, we can finish the proof of
Theorem~\ref{CLR} in the case $r\ne 0$ following the standard argument
of~\cite[Chapters~5.3 and~5.4]{IH81}. More precisely, the weak convergence in
$\Dd_0(\RR)$ of the processes $Z_{n,\vartheta}$ to the process~$Z_\vartheta$
follows from Theorem~5.4.2 of~\cite{Kut84}, which is, in fact, contained
in~\cite{IH81} (without being formulated there).  Note, that the conditions of
this theorem are nothing but Lemmas~\ref{L2}, \ref{L4}
and~\ref{L6}, and that its proof consist in verifying the
condition~\eqref{Skoro}.

\paragraph{Final argument of the proof of Theorem~\ref{CLR} in the case $r=0$}
In this case, it is not possible to establish a lemma similar to
Lemma~\ref{L6}.  In particular, the inequalities~\eqref{dis_1}
and~\eqref{dis_2} do not hold, since in this case (in contrary to the case
$r=0$) the jumps are not becoming seldom.  More precisely, as $n\to+\infty$,
instead of having (on any finite interval) few ``non-vanishing'' jumps, one
has more and more jumps which at the same time become smaller and smaller
(which explains that the trajectories of the limiting likelihood ratio process
in this case are continuous but nowhere differentiable functions).  So, in
order to finish the proof of Theorem~\ref{CLR} we use a different technique.

Since the increments of the process $\ln Z_{n,\vartheta}$ are independent, the
convergence of its restrictions (and hence of those of $Z_{n,\vartheta}$) on
finite intervals $[A,B]\subset\RR$ $\bigl($that is, convergence in the
Skorohod space $\Dd\bigl([A,B]\bigr)$ of functions on $[A,B]$ without
discontinuities of the second kind$\bigr)$ follows from Theorem~6.5.5 of
Gihman and Skorohod~\cite{GS}, Lemma~\ref{L2} and the following lemma.

\begin{lemma}
\label{to_0_L5}
Let the conditions\/ \emph{\textbf{C1}~--~\textbf{C4}} be fulfilled\/ with
$r=0$. Then for any $\varepsilon >0$ we have
$$
\lim_{h \to 0 }\ \lim_{n \to
  +\infty}\ \sup_{\left|u_1-u_2\right|<h} \Pb_\vartheta^{(n)}\bigl(\left|\ln
Z_{n,\vartheta}(u_1)-\ln Z_{n,\vartheta}(u_2)\right| >\varepsilon\bigr)=0.
$$
for all\/ $u_1,u_2\in U_n$ and\/ $\vartheta\in\KK$.
\end{lemma}

Let us note, that taking a closer look on the proof of this lemma, one can see
that we have even a stronger result: for any $\varepsilon >0$ we have
$$
\lim_{h \to 0 }\ \lim_{n \to
  +\infty}\ \sup_{\vartheta\in\KK}\ \sup_{\left|u_1-u_2\right|<h}
\Pb_\vartheta^{(n)}\bigl(\left|\ln Z_{n,\vartheta}(u_1)-\ln Z_{n,\vartheta}(u_2)\right|
>\varepsilon\bigr)=0
$$
for all $u_1,u_2\in U_n$, which allow us to conclude that the convergence of
the restrictions of the process $Z_{n,\vartheta}$ on finite intervals
$[A,B]\subset\RR$ to those of the process $Z_\vartheta$ is uniform with
respect to $\vartheta\in\KK$. Note also, that an alternative way to prove this
convergence (instead of using Lemmas~\ref{L2} and~\ref{to_0_L5}) is to study
the characteristics of the processes and apply, for example, Theorem~7.3.4 of
Jacod and Shiryaev~\cite{JS03}. However, in our opinion, the proof given here
gives more insight on the structure of the considered processes.

In order to conclude the proof of Theorem~\ref{CLR} applying the criterion of
week convergence in $\Dd_0(\RR)$ given in Lemma~\ref{CWC}, we need to check
the condition~\eqref{Skoro}. Since we have already established the convergence
of the restrictions on finite intervals $[A,B]\subset\RR$, it remains to
control the second term of the modulus of continuity
$\Delta_h(Z_{n,\vartheta})$ $\bigl($see, for example,~\cite[Chapters~5.3
  and~5.4]{IH81}$\bigr)$.  So, the last ingredient of the proof of
Theorem~\ref{CLR} is the following estimate on the tails of the process
$Z_{n,\vartheta}$.

\begin{lemma}
\label{to_0_L6}
Let the conditions\/ \emph{\textbf{C1}~--~\textbf{C4}} be fulfilled with\/
$r=0$. Then there exist some constants\/ $b,C>0$ such that
\begin{equation}
\label{to_0_ineqL6}
\Pb_\vartheta^{(n)}\biggl(\sup_{\left|u\right|>D} Z_{n,\vartheta}(u) >e^{-b
  D}\biggr)\leq C e^{-b D}
\end{equation}
for all\/ $D\geq 0$, $n\in\NN$ and\/ $\vartheta\in\KK$.
\end{lemma}

\section{Parameter estimation}
\label{Sec3}

In this section we apply the convergence of normalized likelihood ratio
obtained in Section~\ref{Sec2} to study the problem of parameter estimation
for our model of observations.  In the case $r\ne 0$, the limiting likelihood
ratio being the same as in the fixed jump size case, the properties of
estimators are also the same $\bigl($see, for example,~\cite{Kut84,Kut98} for
more details$\bigr)$.  So, here we consider the case $r=0$ only.

Recall that as function of $\vartheta$, the likelihood of our model given
by~\eqref{LR} is discontinuous (has jumps).  So, the maximum likelihood
estimator $\widehat\vartheta_n$ of $\vartheta$ is introduced through the
equation
$$
\max\Bigl\{L_n\bigl(\,\widehat\vartheta_n+,X^{(n)}\bigr),
L_n\bigl(\,\widehat\vartheta_n-,X^{(n)}\bigr)\Bigr\}=
\sup_{\vartheta\in\Theta}L_n\bigl(\vartheta,X^{(n)}\bigr).
$$

The Bayesian estimator $\widetilde\vartheta_n$ of $\vartheta$ for a given
prior density $p$ and for square loss is defined by
$$
\widetilde\vartheta_n=\frac{\int_{\alpha}^{\beta}\vartheta\,p(\vartheta)
  L_n\bigl(\vartheta,X^{(n)}\bigr)\,\dd \vartheta} {\int_{\alpha}^{\beta} p(\vartheta)
  L_n\bigl(\vartheta,X^{(n)}\bigr)\,\dd \vartheta}\,.
$$

We are interested in the asymptotic properties of the maximum likelihood and
Bayesian estimators of $\vartheta$ as $n\to + \infty$.  To describe the
properties of the estimators we need some additional notations.

We introduce the random variables $\xi_\vartheta$, $\xi^*$, $\zeta_\vartheta$
and $\zeta^*$ by the equations
\begin{align*}
Z_\vartheta(\xi_\vartheta)&=\sup_{u\in\RR}Z_\vartheta(u),\\
Z^*(\xi^*)&=\sup_{u\in\RR}Z^*(u),
\end{align*}
\begin{align*}
\zeta_\vartheta&=\frac{\int_{-\infty}^{+\infty}u\,Z_\vartheta(u)\,\dd u}
           {\int_{-\infty}^{+\infty}Z_\vartheta(u)\,\dd u}\\
\intertext{and}
\zeta^*&=\frac{\int_{-\infty}^{+\infty}u\,Z^*(u)\,\dd u}
     {\int_{-\infty}^{+\infty}Z^*(u)\,\dd u}\:.
\end{align*}

Let us note that $\xi_\vartheta\eqd\psi(\vartheta)\,\xi^*$ and
$\zeta_\vartheta\eqd\psi(\vartheta)\,\zeta^*$.

Now we can state the following theorem giving an asymptotic lower bound on the
risk of all the estimators of $\vartheta$.

\begin{theorem}
\label{to_0_borne}
Let the conditions\/ \emph{\textbf{C1}~--~\textbf{C4}} be fulfilled with\/
$r=0$. Then, for any\/ $\vartheta_0\in\Theta$, we have
$$
\lim_{\delta\to0}\ \Liminf_{n\to+\infty}\ \inf_{\overline
  \vartheta_n}\ \sup_{\left|\vartheta-\vartheta_0\right|<\delta}
\varphi_n^{-2}\,\Ex^{(n)}_\vartheta(\overline\vartheta_n-\vartheta)^2
\geq\Ex\zeta_{\vartheta_0}^2=\psi^2(\vartheta_0)\,\Ex(\zeta^*)^2,
$$
where the\/ $\inf$ is taken over all possible estimators\/ $\overline\vartheta_n$
of the parameter\/ $\vartheta$.
\end{theorem}

This theorem allows us to introduce the following definition.

\begin{definition}
Let the conditions\/ \emph{\textbf{C1}~--~\textbf{C4}} be fulfilled with\/
$r=0$. We say that an estimator\/ $\vartheta_n^*$ is asymptotically efficient
if
$$
\lim_{\delta\to0}\ \lim_{n\to+\infty}\ \sup_{\left|\vartheta -
  \vartheta_0 \right|< \delta}
\varphi_n^{-2}\,\Ex^{(n)}_\vartheta(\vartheta_n^*-\vartheta)^2
=\Ex\zeta_{\vartheta_0}^2=\psi^2(\vartheta_0)\,\Ex(\zeta^*)^2
$$
for all\/ $\vartheta_0\in\Theta$.
\end{definition}

Now, we can state the following two theorems giving the asymptotic properties
of the maximum likelihood and Bayesian estimators.

\begin{theorem}
\label{to_0_Th2}
Let the conditions\/ \emph{\textbf{C1}~--~\textbf{C4}} be fulfilled with\/
$r=0$. Then the maximum likelihood estimator\/ $\widehat{\vartheta}_n$
satisfies uniformly on\/ $\vartheta\in\KK$ the relations
\begin{align*}
\Pb^{(n)}_{\vartheta}-\lim_{n\to+\infty}\widehat\vartheta_n&=\vartheta,\\
\Ll^{(n)}_{\vartheta}\bigl\{\varphi_n^{-1}(\widehat{\vartheta}_n-\vartheta)\bigr\}&
\Rightarrow\Ll(\xi_\vartheta)=\Ll\bigl(\psi(\vartheta)\xi^*\bigr)\\
\intertext{and}
\lim_{n\to+\infty}\Ex^{(n)}_{\vartheta}\varphi_n^{-p}\bigl|\widehat\vartheta_n
-\vartheta\bigr|^p&=\Ex\left|\xi_\vartheta\right|^p=
\psi^p(\vartheta)\,\Ex\left|\xi^*\right|^p \quad\text{for any}\quad p>0.
\end{align*}
In particular, the relative asymptotic efficiency of\/ $\widehat{\vartheta}_n$
is\/ $\Ex(\zeta^*)^2/\Ex(\xi^*)^2$.
\end{theorem}

\begin{theorem}
\label{to_0_Th1}
Let the conditions\/ \emph{\textbf{C1}~--~\textbf{C4}} be fulfilled with\/
$r=0$. Then, for any continuous strictly positive density, the Bayesian
estimator\/ $\widetilde{\vartheta}_n$ satisfies uniformly on\/
$\vartheta\in\KK$ the relations
\begin{align*}
\Pb^{(n)}_{\vartheta}-\lim\limits_{n\to+\infty}\widetilde\vartheta_n
&=\vartheta,\\
\Ll^{(n)}_{\vartheta}\bigl\{\varphi_n^{-1}(\widetilde{\vartheta}_n-\vartheta)\bigr\}&
\Rightarrow\Ll(\zeta_\vartheta)=\Ll\bigl(\psi(\vartheta)\zeta^*\bigr)\\
\intertext{and}
\lim_{n\to+\infty}\Ex^{(n)}_{\vartheta}\varphi_n^{-p}
\bigl|\widetilde\vartheta_n-\vartheta\bigr|^p&=\Ex\left|\zeta_\vartheta\right|^p=
\psi^p(\vartheta)\,\Ex\left|\zeta^*\right|^p\quad\text{for any}\quad p>0.
\end{align*}
In particular,\/ $\widetilde{\vartheta}_n$ is asymptotically efficient.
\end{theorem}

Theorems~\ref{to_0_borne}--\ref{to_0_Th1} follow from the properties of the
normalized likelihood ratio established in Section~\ref{Sec2}. More precisely,
Theorem~\ref{to_0_Th1} is a consequence of Lemmas~\ref{L2}--\ref{L4}
and~\cite[Theorem~1.10.2]{IH81}. Having the properties of the Bayesian
estimators given in Theorem~\ref{to_0_Th1}, we can cite~\cite[Theorem
  1.9.1]{IH81} to provide the proof of Theorem~\ref{to_0_borne}. Finally, the
proof of Theorem~\ref{to_0_Th2} can be carried out following the standard
argument of~\cite[Chapters~5.3 and~5.4]{IH81} which is based on the weak
convergence established in Theorem~\ref{CLR} together with the
inequality~\eqref{to_0_ineqL6}.

\section{Hypothesis testing}
\label{Sec4}

In this section we apply the convergence of normalized likelihood ratio
obtained in Section~\ref{Sec2} to study the problem of hypothesis testing for
our model of observations.  In the case $r\ne 0$, the limiting likelihood
ratio being the same as in the fixed jump size case, the properties of test
are also the same $\bigl($see~\cite{DKY_2} for more details$\bigr)$.  So, here
we consider the case $r=0$ only.

We consider the same model of observation as above, with the only difference
that now we suppose that $\theta\in\Theta=\left[\vartheta_1,b\right)$,
  $0<\vartheta_1<\beta\leq\tau$.  We assume that the conditions
  \textbf{(C1)}--\textbf{(C4)} are fulfilled with $r=0$ and we want to test
  the following two hypothesis:
\begin{align*}
\Hh_1\quad &:\quad \vartheta=\vartheta_1,\\
\Hh_2\quad &:\quad \vartheta>\vartheta_1.
\end{align*}

We define a (randomized) test $\overline\phi _n=\overline\phi
_n\bigl(X^{(n)}\bigr)$ as the probability to accept the hypothesis
$\Hh_2$. The size of the test is defined by
$\Ex_{\vartheta_1}^{(n)}\overline\phi_n\bigl(X^{(n)}\bigr)$, and its power
function is given by
$\beta(\overline\phi_n,\vartheta)=\Ex_\vartheta^{(n)}\overline\phi_n\bigl(X^{(n)}\bigr)$,
$\vartheta>\vartheta_1$. As usually, we denote $\Kk_\varepsilon$ the class of
tests of asymptotic size $\varepsilon\in\left[0,1\right]$, that is,
$$
\Kk_\varepsilon =\Bigl\{\overline\phi_n\quad :\quad \lim_{n\to+\infty}
\Ex_{\vartheta_1}^{(n)}\overline\phi_n\bigl(X^{(n)}\bigr)=\varepsilon\Bigr\}.
$$

Our goal is to construct some tests belonging to this class and to compare
them. The comparison of tests can be done by comparison of their power
functions. It is known that for any reasonable test and for any fixed
alternative the power function tends to $1$. To avoid this difficulty, we use
Pitman's approach and consider \textit{contiguous\/} (or \textit{close}\/)
alternatives.  More precisely, changing the variable
$\vartheta=\vartheta_u\eqdef\vartheta_1+u\varphi_n^*$, where
$\varphi^*_n=\frac{\psi(\vartheta_1)}{n\,r_n^2}$, the initial problem of
hypotheses testing can be replaced by the following one
\begin{align*}
\Hh_1\quad &:\quad u=0,\\
\Hh_2\quad &:\quad u>0,
\end{align*}
and the power function is now
$\beta(\overline\phi_n,u)=\Ex_{\vartheta_n}^{(n)}\overline\phi_n\bigl(X^{(n)}\bigr)$,
$u>0$.

The study is essentially based on the properties of the normalized likelihood
ratio established above.  Note that the limit of the normalized likelihood
ratio at the point $\vartheta=\vartheta_1$ (under hypothesis $\Hh_1$) is the
following:
$$
Z_{n,\vartheta_1}^*(v)=\frac{L_n\bigl(\vartheta_1+v\varphi_n^*,X^{(n)}\bigr)}
{L_n\bigl(\vartheta_1,X^{(n)}\bigr)} \Rightarrow Z^*(v),\quad v\geq 0,
$$
where the process $Z^*$ is defined by~\eqref{Zstar}.

Under alternatives, we obtain
\begin{align*}
Z_{n,\vartheta_1}^*(v)&=\frac{L_n\bigl(\vartheta_1+v\varphi_n^*,X^{(n)}\bigr)}
{L_n\bigl(\vartheta_1,X^{(n)}\bigr)}\\
&=\left(\frac{L_n\bigl(\vartheta_1,X^{(n)}\bigr)}
{L_n\bigl(\vartheta_u,X^{(n)}\bigr)}\right)^{-1}\;
\frac{L_n\bigl(\vartheta_1+v\varphi_n^*,X^{(n)}\bigr)}
     {L_n\bigl(\vartheta_u,X^{(n)}\bigr)}\\
&=\left(\frac{L_n\bigl(\vartheta_u-u\varphi_n^*,X^{(n)}\bigr)}
     {L_n\bigl(\vartheta_u,X^{(n)}\bigr)}\right)^{-1}\;
     \frac{L_n\bigl(\vartheta_u+(v-u)\varphi_n^*,X^{(n)}\bigr)}
          {L_n\bigl(\vartheta_u,X^{(n)}\bigr)}\\
&\Rightarrow \bigl(Z^*(-u)\bigr)^{-1}\,Z^*(v-u)\eqd\exp
          \left\{W(v)-\frac{\left|v-u\right|}{2}+\frac{u}{2}\right\} \eqdef
          Z^*_u(v).
\end{align*}

The score-function test --- which is locally asymptotically uniformly most
powerful (LAUMP) in the regular case $\bigl($see~\cite{DKY_1}$\bigr)$--- does
not exist in this non-regular situation. So, we will construct and study the
general likelihood ratio test (GLRT), Wald's test (WT) and two Bayesian tests
(BT1 and~BT2).

\paragraph{General likelihood ratio test}
The GLRT is defined by the relations
$$
\widehat\phi_n\bigl(X^{(n)}\bigr)=\1_{\bigl\{Q(X^{(n)})>
  h_\varepsilon\bigr\}},
$$
with
$$
Q\bigl(X^{(n)}\bigr)=\sup_{\vartheta >\vartheta_1}\,
\frac{L_n\bigl(\vartheta,X^{(n)}\bigr)} {L_n\bigl(\vartheta_1,X^{(n)}\bigr)}
=\max\left\{\frac{L_n\bigl(\widehat\vartheta_n+,X^{(n)}\bigr)}
{L_n\bigl(\vartheta_1,X^{(n)}\bigr)}\;,\;
\frac{L_n\bigl(\widehat\vartheta_n-,X^{(n)}\bigr)}
     {L_n\bigl(\vartheta_1,X^{(n)}\bigr)}\right\},
$$
where $\widehat\vartheta_n$ is the maximum likelihood estimator of $\theta$.

To choose the threshold $h_\varepsilon$ such that
$\widehat\phi_n\bigl(X^{(n)}\bigr)\in\Kk_\varepsilon$ we need to solve the
following equation (under hypothesis $\Hh_1$)
$$
\Pb_{\vartheta _1}^{\left(n\right)}\left\{Q\bigl(X^{(n)}\bigr)>h_\varepsilon
\right\}= \Pb_{\vartheta _1}^{\left(n\right)}\left\{\sup_{v>0}
Z_{n,\vartheta_1}^*\left(v\right)>h_\varepsilon \right\} \to \Pb
\left\{\sup_{v> 0} Z^*\left(v\right)> h_\varepsilon \right\}=\varepsilon.
$$
For this, we note that the random variable $\sup\limits_{v> 0}\ln
Z^*\left(v\right)$ has the exponential distribution with parameter $1$
$\bigl($see, for example,~\cite{BB95}$\bigr)$. This allows us to calculate
explicitly the threshold $h_\varepsilon $ of the GLRT as solution of the
equation $ 1-e^{-\ln h_\varepsilon}=1-\varepsilon $, that is,
$h_\varepsilon=1/\varepsilon$.

The power function of the GLRT has the following limit:
$$
\beta\bigl(\widehat\phi_n,u\bigr)=\Pb_{\vartheta
  _u}^{\left(n\right)}\left\{\sup_{v> 0}
Z_{n,\vartheta_1}^*\left(v\right)>h_\varepsilon \right\}\to
\Pb\left\{\sup_{v> 0} Z^*_u\left(v\right)>h_\varepsilon \right\}.
$$
This limiting power function is obtained below with the help of numerical
simulations.

\paragraph{Wald's test}
To define the WT, let us note that the maximum likelihood estimator
$\widehat\vartheta_n$ converges in distribution:
$$
\left(\varphi_n^*\right)^{-1}\bigl(\widehat\vartheta_n-\vartheta
_1\bigr)\Rightarrow \xi_+^*,
$$
where the random variable $\xi^*_+$ is solution of the equation
$$
Z^*(\xi_+^*)=\sup_{v>0}Z^*(v).
$$
Therefore, if we put
$$
\phi_n^\circ\bigl(X^{(n)}\bigr)=\1_{\left\{\left(\varphi_n^*\right)^{-1}\bigl(\widehat\vartheta
  _n-\vartheta _1\bigr) >m_\varepsilon \right\}},
$$
where $m_\varepsilon$ is defined by the equation
$$
\Pb\left\{\xi_+^*>m_\varepsilon\right\}=\varepsilon,
$$
then $\phi_n^\circ\in\Kk_\varepsilon$.

We recall the result of~\cite{Fu07}, that the joint distribution of $\bigl(\ln
Z^*\left(\xi_+^*\right),\xi_+^*\bigr)$ has the density
$$
f(y,t)=\frac{y}{\sqrt{2\pi
    t^3}}\exp\left\{-\frac{\left(y+\frac{t}{2}\right)^2}{2t}\right\},
$$
which allows us to calculate the marginal density of $\xi_+^*$ as follows:
\begin{align*}
f(t)&=\int_0^{+\infty}f(y,t)\,\dd
y=\int_0^{+\infty}\frac{\frac{y}{\sqrt{t}}}{\sqrt{2\pi t}}\,\exp
\left\{-\frac12\left(\frac{y}{\sqrt{t}}+\frac{\sqrt{t}}{2}\right)^2\right\}
\;\dd\left(\frac{y}{\sqrt{t}}\right)\\
&=\int_0^{+\infty}\frac{z}{\sqrt{2\pi t}}\,\exp
\left\{-\frac12\left(z+\frac{\sqrt{t}}{2}\right)^2\right\}\,\dd z
=\int_{\frac{\sqrt{t}}{2}}^{+\infty}\frac{x-\frac{\sqrt{t}}{2}}{\sqrt{2\pi
    t}}\,\exp\left\{-\frac{x^2}2\right\}\,\dd x\\
&=-\int_{\frac{\sqrt{t}}{2}}^{+\infty}\frac1{\sqrt{2\pi
    t}}\;\dd\exp\left\{-\frac{x^2}2\right\}
-\int_{\frac{\sqrt{t}}{2}}^{+\infty}\frac{\frac{\sqrt{t}}{2}}{\sqrt{2\pi
    t}}\,\exp\left\{-\frac{x^2}2\right\}\,\dd x\\
&=\frac1{\sqrt{2\pi t}}\,\exp\left\{-\frac{t}8\right\}
-\frac12\Phi\left(-\frac{\sqrt{t}}{2}\,\right),
\end{align*}
where $\Phi$ is the distribution function of the standard Gaussian low
$\Nn(0,1)$. So, the threshold $m_\varepsilon$ can be obtained as the solution
of the equation
\begin{equation}
\label{WT_thr_eq}
\int_{m_\varepsilon}^{+\infty}\left(\frac1{\sqrt{2\pi t}}
\exp\left\{-\frac{t}8\right\}
-\frac12\Phi\left(-\frac{\sqrt{t}}{2}\,\right)\right)\,\dd t=\varepsilon.
\end{equation}

The power function of the WT has the following limit:
$$
\beta\bigl(\phi_n^\circ,u\bigr)=\Pb_{\vartheta _u}^{(n)}
\left\{\left(\varphi_n^*\right)^{-1}\bigl(\widehat\vartheta _n-\vartheta
_u\bigr)+u>m_\varepsilon \right\}\to \Pb\left\{\xi_u^*
>m_\varepsilon-u\right\},
$$
where the random variable $\xi_u^*$ is solution of the equation
$$
Z\left(\xi_u^*\right)=\sup_{v> -u } Z^*\left(v\right).
$$

Note that we can also derive another expression of the limiting power function
of the WT as follows:
$$
\beta\bigl(\phi_n^\circ,u\bigr)=\Pb_{\vartheta
  _u}^{(n)}\left\{\left(\varphi_n^*\right)^{-1}\bigl(\widehat\vartheta
_n-\vartheta _1\bigr) >m_\varepsilon
\right\}\to\Pb\left\{\xi_{u,+}^* >m_\varepsilon\right\},
$$
where the random variable $\xi_{u,+}^*$ is solution of the equation
$$
Z\left(\xi_{u,+}^*\right)=\sup_{v>0} Z^*_u\left( v\right).
$$
The threshold and the limiting power function are obtained below with the help
of numerical simulations.

\paragraph{Bayesian tests}
Suppose now that the parameter $\vartheta $ is a random variable with the
\textit{a priori\/} density $p(\theta)$, $\vartheta _1\leq \theta
<\beta$. This density is supposed to be continuous and positive. We consider
two tests.

The first one (BT1) is based on the Bayesian estimator:
$$
\widetilde\phi_n\bigl(X^{(n)}\bigr)=
\1_{\left\{\left(\varphi_n^*\right)^{-1}\left(\widetilde\vartheta _n-\vartheta
  _1\right) >k_\varepsilon \right\}}.
$$
As above, we have the convergence in distribution:
$$
\left(\varphi_n^*\right)^{-1}\bigl(\widetilde\vartheta_n-\vartheta
_1\bigr)\Rightarrow\zeta_+^*\eqdef\frac{\int_{0}^{+\infty }v\,
  Z^*\left(v\right)\,\dd v}{\int_{0}^{+\infty }Z^*\left(v\right)\,\dd v}\,,
$$
which allows us to chose the threshold such that
$\widetilde\phi_n\in\Kk_\varepsilon$ as the solution of the equation
\begin{equation}
\label{BT1_thr_eq}
\Pb\left\{\zeta_+^*>k_\varepsilon \right\}=\varepsilon.
\end{equation}

The power function of the BT1 has the following limit:
$$
\beta\bigl(\widetilde\phi_n,u\bigr)=\Pb_{\vartheta _u}^{(n)}
\left\{\left(\varphi_n^*\right)^{-1}\bigl(\widetilde\vartheta _n-\vartheta
_u\bigr)+u>k_\varepsilon \right\}\to \Pb\left\{\zeta_u^*
>k_\varepsilon-u\right\},
$$
where the random variable $\zeta_u^*$ is given by
$$
\zeta_u^*= \frac{\int_{-u}^{+\infty }v\,Z^*\left(v\right)\,\dd
  v}{\int_{-u}^{+\infty }Z^*\left(v\right)\,\dd v}.
$$

Note that we can also derive another expression of the limiting power function
of the BT1 as follows:
$$
\beta\bigl(\widetilde\phi_n,u\bigr)=\Pb_{\vartheta
  _u}^{(n)}\left\{\left(\varphi_n^*\right)^{-1}\bigl(\widehat\vartheta
_n-\vartheta _1\bigr) > k_\varepsilon
\right\}\to\Pb\left\{\zeta_{u,+}^* >k_\varepsilon\right\},
$$
where the random variable $\zeta_{u,+}^*$ is given by
$$
\zeta_{u,+}^*=\frac{\int_{0}^{+\infty }v\,Z^*_u\left(v\right)\,\dd
  v}{\int_{0}^{+\infty }Z^*_u\left(v\right)\,\dd v}
$$
The threshold and the limiting power function are obtained below with the help
of numerical simulations.

The second test (BT2) minimizes the mean error. The likelihood ratio is
$$
\widetilde L\bigl(X^{(n)}\bigr)=\int_{\vartheta
  _1}^{\beta}\frac{L_n\bigl(\theta,X^{(n)}\bigr)}{L_n\bigl(\vartheta
  _1,X^{(n)}\bigr)}\,p(\theta)\,\dd\theta
=\varphi_n^*\int_{0}^{\left(\varphi_n^*\right)^{-1}\,(\beta-\vartheta_1)}
Z^*_{n,\vartheta_1}(v)\,p(\vartheta _1+v\varphi_n^*)\,\dd v.
$$
Hence, we have the following limit:
$$
(\varphi_n^*)^{-1}\,\widetilde L\bigl(X^{(n)}\bigr)\Rightarrow
p\left(\vartheta _1\right)\int_{0}^{+\infty } \exp\left\{
W\left(v\right)-\frac{v}{2} \right\}\,\dd v.
$$
Therefore, if we denote
$$
R_n=\frac{(\varphi_n^*)^{-1}\, \widetilde
  L\bigl(X^{(n)}\bigr)}{p\left(\vartheta _1\right)}
$$
and chose $g_\varepsilon$ as solution of the equation
$$
\Pb\left\{ \int_{0}^{+\infty} \exp\left\{ W\left(v\right)-\frac{v}{2}
\right\}\,\dd v >g_\varepsilon \right\}=\varepsilon,
$$
the test $\1_{\left\{R_n>g_\varepsilon \right\}}$ belongs to the class
$\Kk_\varepsilon $.

\begin{figure}[ht]
\centerline{\includegraphics[width=\textwidth]{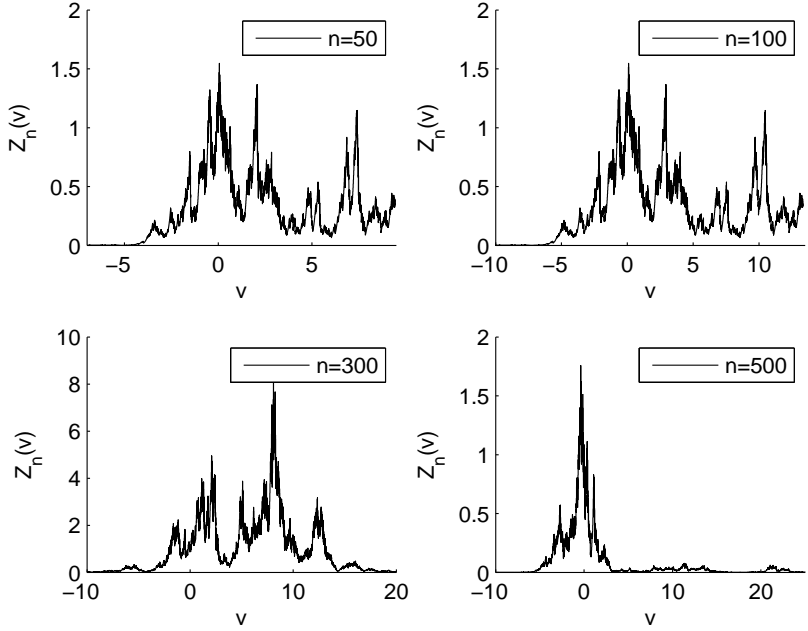}}
\caption{Some realization of $Z_{n,\vartheta_1}^*(v)$}
\label{Z_u_n_cp_0}
\end{figure}

\paragraph{Numerical simulations}
Now, let us carry out some numerical simulations for the GLRT, the WT and the
BT1.  We take $r_n=n^{-0.25}$ and, in order to simplify the simulations, we
take a function $\psi_n(t)$ depending neither on $n$ nor on $t$.  More
precisely, we consider $n$ independent trajectories
$X^{(n)}_j=\bigl\{X^{(n)}_j(t),\ t\in\left[0,4\right]\bigr\}$, $j=1,\ldots,n$,
of an inhomogeneous Poisson process on the interval $\left[0,4\right]$ of
intensity function
$$
\lambda_{\vartheta}^{(n)}(t)=1.5+n^{-0.25}\,\1_{\left\{t>
  \vartheta\right\}},\qquad 0\leq t\leq 4,
$$
with $\vartheta \in \left[2,4\right)$.  So, denoting $\vartheta_1=2$ and
$$
\varphi_n^*=\frac{\psi(\vartheta_1)}{n\,r_n^2}=\frac{1.5}{\sqrt{n}}\,,
$$
we have (for $v\geq 0$)
\begin{align*}
\ln Z_{n,\vartheta_1}^*(v)&=
\sum_{j=1}^n\int_{\left(\vartheta_1,\vartheta_1+v\varphi_n^*\right]} \ln
\frac{1.5}{1.5+n^{-0.25}}\, X_j^{(n)}(\dd t)+1.5\,v\,n^{0.25}\\
&=\ln \frac{1.5}{1.5+n^{-0.25}}
\sum_{j=1}^n\,\left(X_j^{(n)}(\vartheta_1+v\varphi_n^*)
-X_j^{(n)}(\vartheta_1)\right) +1.5\,v\, n^{0.25}.
\end{align*}
Some realizations of $Z_{n,\vartheta_1}^*$ can be found in Figure~\ref{Z_u_n_cp_0}.

Recall that the threshold $h_\varepsilon=1/\varepsilon$ of the GLRT is known
explicitly. We obtain the threshold $m_\varepsilon$ of the WT by numerically
solving the equation~\eqref{WT_thr_eq}, while the threshold $k_\varepsilon$ of
the BT1 is obtained from the equation~\eqref{BT1_thr_eq} by means of numerical
simulations of the random variable~$\zeta^*_+$. Some values of the thresholds
$m_\varepsilon$ and $k_\varepsilon$ are given in Table~\ref{Thr_cp_0}.

\begin{table}[ht]
\caption{Thresholds of WT and BT1}
\label{Thr_cp_0}
\centerline{%
\begin{tabular}{|c|c|c|c|c|c|c|}
 \hline
 $\varepsilon$ &  0.001 &  0.005 &  0.01  & 0.05  & 0.1   & 0.2   \\
 \hline
 $m_\varepsilon$  & 30.336 & 20.686 & 14.886 & 7.282 & 4.531 & 2.236 \\
 \hline
 $k_\varepsilon$  & 24.877 & 17.588 & 16.782 & 8.582 & 5.573 & 3.024 \\
 \hline
\end{tabular}}
\end{table}

\begin{figure}[ht]
\centerline{\includegraphics[width=\textwidth]{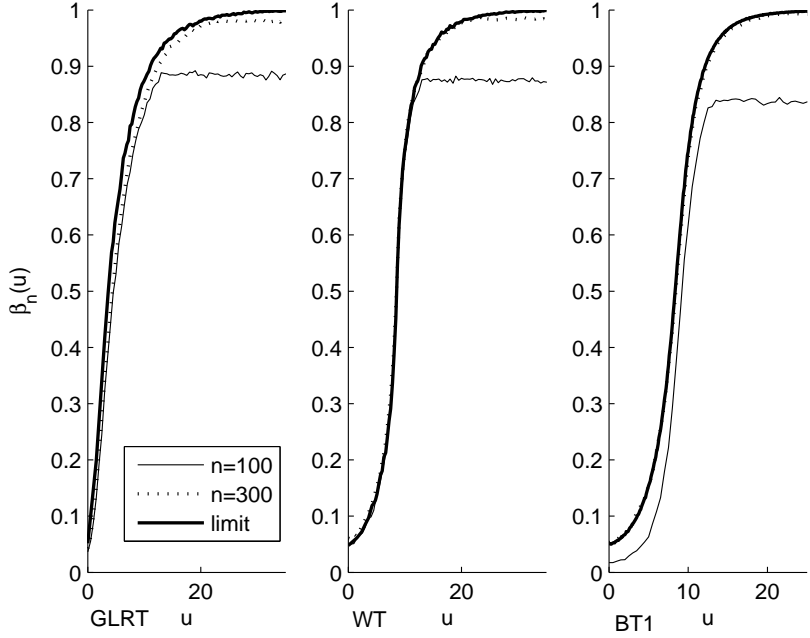}}
\caption{Power functions of GLRT, WT and BT1}
\label{PF_cp_0_n}
\end{figure}

To illustrate the convergence of power functions of different tests to their
limits, we present in Figure~\ref{PF_cp_0_n} the power functions for $n=100$
($r_n=0.3162$) and $n=300$ ($r_n=0.2403$), as well as the limiting power
functions.  All these power functions are obtained by means of numerical
simulations.  Note that the values of $u$ greater than $2(\varphi_n^*)^{-1}$
correspond to $\theta_u=\theta_1+u\varphi_n^*>4$, which means that there is no
longer jump in intensity function on the interval $\left[0,4\right]$.  This
explains the fact that for $n=100$, the power functions are constant for
$u>2(\varphi_{100}^*)^{-1}\approx 13.33$.

\paragraph{Comparison of the limiting power functions}
Let us fix an alternative $u_1>0$ and consider the testing problem with two
simple hypotheses
\begin{align*}
\Hh_1^{\phantom{u_1}}\quad &:\qquad u=0,\\
\Hh_2^{u_1}\quad &:\qquad u=u_1.
\end{align*}
Remind that in this situation the most powerful test is the Neyman-Pearson
test (N-PT).  Of course, it is impossible to use the N-PT in our initial
problem, because it depends on the value $u_1$ under alternative which is
unknown.  However, its power (considered as function of $u_1$) gives an upper
bound (Neyman-Pearson envelope) for the power functions of all the tests.
Therefore, it is interesting to compare the power functions of different tests
not only one with another, but also with the power of the N-PT.

The N-PT is given by
$$
\phi_n^*\bigl(X^{(n)}\bigr)=\1_{\bigl\{Z_{n,\vartheta_1}^*(u_1)>
  d_\varepsilon\bigr\}}+q_\varepsilon\1_{\bigl\{Z_{n,\vartheta_1}^*(u_1)=
  d_\varepsilon\bigr\}},
$$
where $d_\varepsilon$ and $q_\varepsilon$ are solution of the equation
\begin{equation}
\label{NPT_eq}
\Pb_{\vartheta_1}^{(n)}\bigl(Z_{n,\vartheta_1}^*(u_1)>d_\varepsilon\bigr)+q_\varepsilon\,
\Pb_{\vartheta_1}^{(n)}\bigl(Z_{n,\vartheta_1}^*(u_1)=d_\varepsilon\bigr)=\varepsilon.
\end{equation}

Recall that the likelihood ratio $Z_{n,\vartheta_1}^*(u_1)$ under hypothesis
$\Hh_1$ converges to the following limit
$$
Z_{n,\vartheta_1}^*(u_1) \Rightarrow Z^*(u_1)
=\exp\left\{W\left(u_1\right)-\frac{u_1}{2}\right\}.
$$
Hence, in the asymptotic setting, the equation~\eqref{NPT_eq} can be replaced
by the equation
$$
\Pb\bigl(Z^*(u_1)>d_\varepsilon\bigr)+q_\varepsilon\,
\Pb\bigl(Z^*(u_1)=d_\varepsilon\bigr)=\varepsilon
$$
and, since $Z^*(u_1)$ is a continuous random variable, we can put
$q_\varepsilon = 0$ and find the threshold $d_\varepsilon $ as the solution of
the equation
$$
\Pb\bigl(Z^*(u_1)> d_\varepsilon\bigr)=\varepsilon.
$$
Note that
$$
\Pb\bigl(Z^*(u_1)> d_\varepsilon\bigr)=\Pb\left(W(u_1)> \ln
d_\varepsilon+\frac{u_1}2\right)=\Pb\left(\zeta>\frac{\ln
  d_\varepsilon+\frac{u_1}2}{\sqrt{u_1}}\right),
$$
where $\zeta \sim \Nn\left(0,1\right)$. Therefore, denoting $z_\varepsilon$
the quantile of order $1-\varepsilon$ of the standard Gaussian law
$\bigl(\Pb\left(\zeta >z_\varepsilon \right)=\varepsilon\bigr)$, the threshold
$d_\varepsilon$ is given by
$$
d_\varepsilon =e^{z_\varepsilon \sqrt{u_1}-\frac{u_1}{2}}.
$$

Under alternative $\Hh_2^{u_1}$, we have
$$
Z_{n,\vartheta_1}^*(u_1)=\frac{L_n\bigl(\vartheta_1+u_1\varphi_n^*,X^{(n)}\bigr)}
{L_n\bigl(\vartheta_1,X^{(n)}\bigr)}
=\biggl(\frac{L_n\bigl(\vartheta_1+u_1\varphi_n^*-u_1\varphi_n^*,X^{(n)}\bigr)}
{L_n\bigl(\vartheta_1+u_1\varphi_n^*,X^{(n)}\bigr)}\biggr)^{-1}\Rightarrow
\bigl(Z^*(-u_1)\bigr)^{-1},
$$
which allows us to obtain the limiting power of the N-PT as follows:
\begin{align*}
\beta\bigl(\phi_n^*\bigr)&=\Pb_{\vartheta_1+u_1\varphi_n^*}^{(n)}
\left(Z_{n,\vartheta_1}^*(u_1)>d_\varepsilon\right)\\
&\to \Pb\left(\bigl(Z^*(-u_1)\bigr)^{-1}>d_\varepsilon\right)
=\Pb\left(\exp\left\{-W\left(-u_1\right)+\frac{u_1}2\right\}>d_\varepsilon\right)\\
&=\Pb\left(W\left(u_1\right)>\ln
d_\varepsilon-\frac{u_1}2\right)=\Pb\left(\zeta>\frac{\ln
  d_\varepsilon-\frac{u_1}2}{\sqrt{u_1}}\right)=\Pb\bigl(\zeta
>z_\varepsilon-\sqrt{u_1}\bigr).
\end{align*}
So, the limiting Neyman-Pearson envelope is given by
$$
\beta(u)=\Pb\bigl(\zeta >z_\varepsilon-\sqrt{u}\bigr)=
1-\Phi\bigl(z_\varepsilon-\sqrt{u}\bigr),
$$
where, as before, $\Phi$ is the distribution function of the standard Gaussian
low.

\begin{figure}[ht]
\centerline{\includegraphics[width=\textwidth]{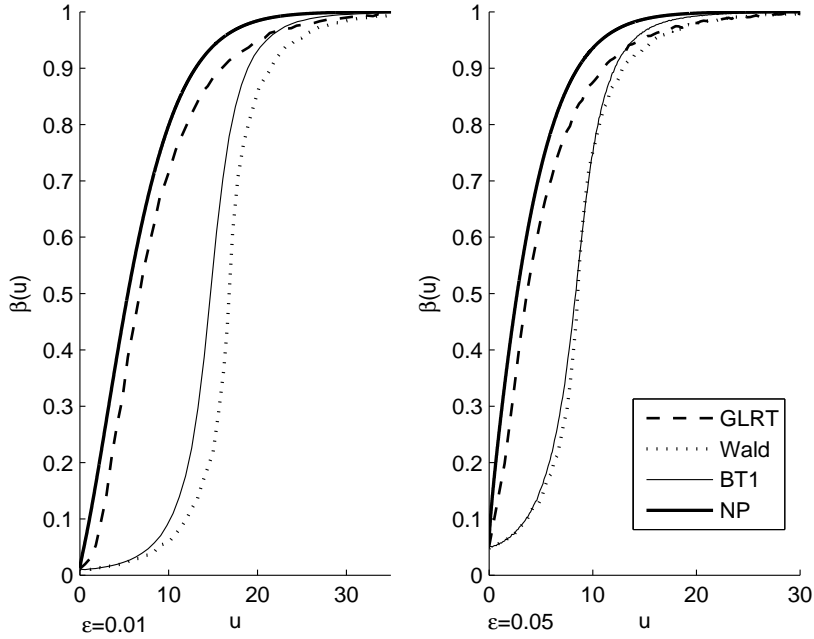}}
\caption{Comparison of limiting power functions for $\varepsilon=0.05$ and
  $\varepsilon=0.4$}
\label{PF_cp_0_comp}
\end{figure}

The limiting power functions of the GLRT, of the WT and of the BT1 are
obtained by means of numerical simulations and are presented in
Figure~\ref{PF_cp_0_comp} together with the limiting Neyman-Pearson envelope
$\beta(u)$.

We can observe that the limiting power function of the GLRT is the closest to
the limiting Neyman-Pearson envelope for small values of $u$, while the
limiting power function of the BT1 is the one that tends to $1$ (as $u$
becomes large) the most quickly.  We can also see that for $\varepsilon=0.05$
the limiting power functions of the WT and of the BT1 are close (especially
when $u$ is small).  Finally, we need to say that all these limiting power
functions are perceptibly below the limiting Neyman-Pearson envelope, and that
the choice of the asymptotically optimal test remains an open question.

\section{Proofs of the lemmas}
\label{Sec5}

The proofs of Lemmas~\ref{L2}--\ref{L4} in the case $r\ne 0$, as
well as the proof of Lemma~\ref{L6}, are similar to the fixed jump size
case and hence are omitted $\bigl($the interested reader can see, for
example,~\cite{Kut84,Kut98}$\bigr)$.

\begin{lemproof}[ in the case $r=0$]{L2}
First we study the convergence of 2-dimensional distributions. For this,
consider the distribution of the vector
$\bigl(Z_{n,\vartheta}(u_1),Z_{n,\vartheta}(u_2)\bigr)$ with some fixed
$u_1,u_2\in\RR$. The characteristic function of the natural logarithm of this
vector can be written as follows $\bigl($see, for
example,~\cite{Kut84}$\bigr)$:
\begin{align*}
&\Ex^{(n)}_{\vartheta}\exp\bigl(it_1\ln Z_{n,\vartheta}(u_1)+it_2\ln
  Z_{n,\vartheta}(u_2)\bigr)\\
&\qquad\qquad=\exp\biggl\{n\int_0^{\tau}\biggl(\exp\Bigl\{it_1\ln
  \frac{\lambda_{\vartheta+u_1\varphi_n}^{(n)}(t)}{\lambda_\vartheta^{(n)}(t)}+it_2\ln
  \frac{\lambda_{\vartheta+u_2\varphi_n}^{(n)}(t)}{\lambda_\vartheta^{(n)}(t)}\Bigr\}-1
  \\
&\ \qquad\qquad\qquad\qquad\qquad
  -it_1\Bigl(\frac{\lambda_{\vartheta+u_1\varphi_n}^{(n)}(t)}
  {\lambda_\vartheta^{(n)}(t)}-1\Bigr)
  -it_2\Bigl(\frac{\lambda_{\vartheta+u_2\varphi_n}^{(n)}(t)}{\lambda_\vartheta^{(n)}(t)}-1\Bigr)\biggr)
  \lambda_\vartheta^{(n)}(t)\,\dd t\biggr\}\\
&\qquad\qquad=\exp\Bigl\{A_{n,\vartheta}(u_1,u_2,t)\,\Bigr\}
\end{align*}
with an evident notation.

We will consider the case $u_2> u_1\geq 0$ only (the other cases can be treated
in a similar way). In this case, we have
\begin{align*}
A_{n,\vartheta}(u_1,u_2,t)&=n\int_{\vartheta}^{\vartheta+u_1\varphi_n
}\biggl(\exp\Bigl\{(it_1+it_2)\ln \frac{\psi_n(t)}{\psi_n(t)+r_n}\Bigr\}-1\\
&\ \qquad\qquad\qquad-(it_1+it_2)\Bigl(\frac{\psi_n(t)}{\psi_n(t)+r_n}-1\Bigr)\biggr)
\Bigl(\psi_n(t)+r_n\Bigr)\,\dd t \\
&\ \quad+n\int_{\vartheta+u_1\varphi_n }^{\vartheta+u_2\varphi_n
}\biggl(\exp\Bigl\{it_2\ln \frac{\psi_n(t)}{\psi_n(t)+r_n}\Bigr\}-1\\
&\ \qquad\qquad\qquad\quad-it_2\Bigl(\frac{\psi_n(t)}{\psi_n(t)+r_n}-1\Bigr)\biggr)
\Bigl(\psi_n(t)+r_n\Bigr)\,\dd t\\
&=n I_1+n I_2
\end{align*}
with evident notations.

Using the mean value theorem for the integrals $I_1$ and $I_2$, it is possible
to find some $s_n \in \left(\vartheta,\vartheta+u_1\varphi_n \right)$ and $v_n \in
    \left(\vartheta+ u_1\varphi_n,\vartheta+u_2\varphi_n \right)$ such that
\begin{align*}
n I_1&=\frac{u_1}{r_n^2}\Bigl(\exp\{i(t_1+t_2)\ln g_n(s_n)\}-1
-i(t_1+t_2)\bigl(g_n(s_n)-1\bigr)\Bigr) \Bigl(\psi_n(s_n)+r_n\Bigr)\\
\intertext{and}
n I_2&=\frac{u_2-u_1}{r_n^2}\Bigl(\exp\{it_2\ln g_n(v_n)\}-1
-it_2\bigl(g_n(v_n)-1\bigr)\Bigr) \Bigl(\psi_n(v_n)+r_n\Bigr),
\end{align*}
where we have denoted
$g_n(t)=\frac{\psi_n(t)}{\psi_n(t)+r_n}=1-\frac{r_n}{\psi_n(t)+r_n}\,$.

As $s_n\to\vartheta$, using the condition \textbf{C3} we obtain $\lim\limits_{n
  \to +\infty} \psi_n(s_n)=\psi(\vartheta)$. So,
$$
n I_1\sim\frac{u_1\psi(\vartheta)}{r_n^2}\Bigl(\exp\{i(t_1+t_2)\ln g_n(s_n)\}-1
-i(t_1+t_2)\bigl(g_n(s_n)-1\bigr)\Bigr).
$$

As $r_n \to 0$ and $\ell\leq\psi _n(t)+r_n\leq L$, we have
$g_n(s_n)-1=O(r_n)\to0$. So, using Taylor expansion we get
\begin{align*}
\ln g_n(s_n)&=\ln\Bigr(1+\bigl(g_n(s_n)-1\bigr)\Bigr)\\
&=g_n(s_n)-1-\frac12\bigl(g_n(s_n)-1\bigr)^2
+o\biggl(\frac{r_n^2}{(\psi_n(s_n)+r_n)^2}\biggr)\\
&=g_n(s_n)-1-\frac12\bigl(g_n(s_n)-1\bigr)^2+o(r_n^2).
\end{align*}
In particular, $\ln g_n(s_n)=O(r_n)$ and $\bigl(\ln
g_n(s_n)\bigr)^2=\bigl(g_n(s_n)-1\bigr)^2+o(r_n^2)$.

Using Taylor expansion once more, we obtain
$$
\exp\bigl(it\ln g_n(s_n)\bigr)=1+it\ln g_n(s_n)-\frac{t^2}{2}\bigl(\ln
g_n(s_n)\bigr)^2+o(r_n^2).
$$
So,
\begin{align*}
n I_1&\sim\frac{u_1\psi(\vartheta)}{r_n^2}\Bigl(-i(t_1+t_2)
\frac{\bigl(g_n(s_n)-1\bigr)^2}{2}-
\frac{(t_1+t_2)^2}{2}\bigl(g_n(s_n)-1\bigr)^2+o(r_n^2)\Bigr)\\
&=\frac{u_1\psi(\vartheta)}{r_n^2}\Bigl(-\frac{i(t_1+t_2)r_n^2}{2(\psi(s_n)+r_n)^2}-
\frac{(t_1+t_2)^2\,r_n^2} {2(\psi(s_n)+r_n)^2}+o(r_n^2)\Bigr)\\
& \to
\frac{u_1}{\psi(\vartheta)}\Bigl(-\frac{i(t_1+t_2)}{2}-\frac{(t_1+t_2)^2}{2}\Bigr).
\end{align*}

Similarly, we can show that
$$
n I_2 \to
\frac{u_2-u_1}{\psi(\vartheta)}\Bigl(-\frac{it_2}{2}-\frac{t_2^2}{2}\Bigr),
$$
and hence
\begin{equation}
\label{to_0_chfunc1}
\begin{aligned}
&\Ex^{(n)}_{\vartheta}\exp\bigl(it_1\ln Z_{n,\vartheta}(u_1)+it_2\ln Z_{n,\vartheta}(u_2)\bigr)\\
&\qquad\qquad\to
\exp\Bigl\{-\frac{u_2-u_1}{2\psi(\vartheta)}(it_2+t_2^2)-
\frac{u_1}{2\psi(\vartheta)}\bigl(i(t_1+t_2)+(t_1+t_2)^2\bigr)\Bigr\}.
\end{aligned}
\end{equation}

For all $u>0$, we introduce the $\sigma$-algebra
$\Ff_{u}=\sigma\bigl\{W(v),\ 0\leq v\leq u\bigr\}$ and write
\begin{align*}
&\Ex\exp\bigl(it_1\ln Z_\vartheta(u_1)+it_2\ln Z_\vartheta(u_2)\bigr)\\
&\qquad\qquad=\Ex\biggl(\exp\bigl\{i(t_1+t_2)\ln Z_\vartheta(u_1)\bigr\}
  \Ex\Bigl(\exp\bigl\{it_2\bigl(\ln Z_\vartheta(u_2)-\ln
  Z_\vartheta(u_1)\bigr)\bigr\}\Bigm|\Ff_{u_1}\Bigr)\biggr)\\
&\qquad\qquad=\exp\Bigl\{-\frac{(t_1+t_2)^2}{2\psi(\vartheta)}u_1
  -\frac{i(t_1+t_2)}{2\psi(\vartheta)}u_1-\frac{t_2^2}{2\psi(\vartheta)}(u_2-u_1)
  -\frac{it_2}{2\psi(\vartheta)}(u_2-u_1)\Bigr\}.
\end{align*}
Combining this with~\eqref{to_0_chfunc1}, we obtain the convergence of
2-dimensional distributions.  The convergence of three and more dimensional
distributions can be carried out in a similar way, and the uniformity with
respect to $\vartheta$ is obvious.
\end{lemproof}

\begin{lemproof}[ in the case $r=0$]{L3}
We will consider the case $u_2\geq u_1\geq 0$ only (the other cases can be
treated in a similar way). According to~\cite[Lemma 1.1.5]{Kut98}, we have
\begin{align*}
\Ex^{(n)}_{\vartheta}\bigl|Z_{n,\vartheta}^{1/2}(u_1)-Z_{n,\vartheta}^{1/2}(u_2)\bigr|^2
&\leq n\int_{0}^{\tau}\biggl(\sqrt{\lambda_{\vartheta+u_1\varphi_n}^{(n)}(t)}
  -\sqrt{\lambda_{\vartheta+u_2\varphi_n}^{(n)}(t)}\,\biggr)^2\,\dd t\\
&=n\int_{\vartheta+ u_1\varphi_n }^{\vartheta+ u_2\varphi_n
}\bigl(\sqrt{\psi_n(t)+r_n}-\sqrt{\psi_n(t)}\,\bigr)^2\,\dd t\\
&=n\int_{\vartheta+ u_1\varphi_n }^{\vartheta+ u_2\varphi_n
}\frac{r_n^2}{\bigl(\sqrt{\psi_n(t)+r_n}
  +\sqrt{\psi_n(t)}\,\bigr)^2}\;\dd t.
\end{align*}

As $\lambda_\vartheta^{(n)}$ is uniformly separated from zero, we have
$$
\bigl(\sqrt{\psi_n(t)+r_n} +\sqrt{\psi_n(t)}\,\bigr)^2\geq \bigl(\sqrt{\ell}
+\sqrt{\ell}\,\bigr)^2=4\ell,
$$
and hence
$$
\Ex^{(n)}_{\vartheta}\bigl|Z_{n,\vartheta}^{1/2}(u_1)-Z_{n,\vartheta}^{1/2}(u_2)\bigr|^2
\leq n\int_{\vartheta+ u_1\varphi_n }^{\vartheta+ u_2\varphi_n}\,
\frac{r_n^2}{4\ell}\;\dd t =\frac{1}{4\ell}\left|u_1-u_2\right|.
$$
So, the required inequality holds with $C=\frac{1}{4\ell}$.
\end{lemproof}

\begin{lemproof}[ in the case $r=0$]{L4}
We will consider the case $u\geq 0$ only (the other case can be treated in a
similar way). According to~\cite[Lemma 1.1.5]{Kut98}, we have
\begin{align*}
\Ex^{(n)}_{\vartheta}Z_{n,\vartheta}^{1/2}(u)&=\exp\Biggl\{-\frac{n}{2}\int_0^{\tau}
\biggl(\sqrt{\lambda_{\vartheta+u\varphi_n}^{(n)}(t)}
-\sqrt{\lambda_\vartheta^{(n)}(t)}\,\biggr)^2\,\dd t\Biggr\}\\
&=\exp\Biggl\{-\frac{n}{2}\int_{\vartheta}^{\vartheta+ u\varphi_n}
\bigl(\sqrt{\psi_n(t)}
  -\sqrt{\psi_n(t)+r_n}\,\bigr)^2\,\dd t\Biggr\}\\
&=\exp\biggl\{-\frac{n}{2}\int_{\vartheta}^{\vartheta+ u\varphi_n
}\frac{r_n^2}{\bigl(\sqrt{\psi_n(t)}
  +\sqrt{\psi_n(t)+r_n}\,\bigr)^2}\;\dd t\biggr\}.
\end{align*}

As $\lambda_\vartheta^{(n)}$ is uniformly bounded, we have
$$
\bigl(\sqrt{\psi_n(t)+r_n}
  +\sqrt{\psi_n(t)}\,\bigr)^2
 \leq \bigl(\sqrt{L}
  +\sqrt{L}\,\bigr)^2 =4L,
$$
and hence
$$
\Ex^{(n)}_{\vartheta}Z_{n,\vartheta}^{1/2}(u)
\leq\exp\Bigl\{-\frac{n}{2}\int_{\vartheta}^{\vartheta+ \varphi_n
u}\frac{r_n^2}{4L}\,\dd t\Bigr\}=\exp\Bigl\{-\frac{1}{8L}\left|u\right|\Bigr\}.
$$
So, the required inequality holds with $k_*=\frac{1}{8L}$.
\end{lemproof}

\begin{lemproof}{to_0_L5}
Using Markov inequality, we get
$$
\Pb_\vartheta^{(n)}\bigl(\left|\ln Z_{n,\vartheta}(u_1)-\ln Z_{n,\vartheta}(u_2)\right|
>\varepsilon\bigr) \leq \frac{1}{\varepsilon^2}\,\Ex_\vartheta^{(n)}\bigl(\ln
Z_{n,\vartheta}(u_1)-\ln Z_{n,\vartheta}(u_2)\bigr)^2.
$$

First we consider the case $u_1,u_2\geq 0$ (and say $u_2\geq u_1$). In this
case, we have
\begin{align*}
\ln Z_{n,\vartheta}(u_2)-\ln Z_{n,\vartheta}(u_1)
&=\sum_{j=1}^n\int_{\vartheta+u_1\varphi_n}^{\vartheta+u_2\varphi_n}\ln
\frac{\psi_n(t)}{\psi_n(t)+r_n}\,\dd
X_j^{(n)}(t)+n\int_{\vartheta+u_1\varphi_n}^{\vartheta+u_2\varphi_n}r_n\,\dd
t\\
&=\sum_{j=1}^n
\int_{\vartheta+u_1\varphi_n}^{\vartheta+u_2\varphi_n}\ln\frac{\psi_n(t)}{\psi_n(t)+r_n}\,\dd
Y_j^{(n)}(t)\\
&\ \quad+n\int_{\vartheta+u_1\varphi_n}^{\vartheta+u_2\varphi_n}\Bigl(
\bigl(\psi_n(t)+r_n\bigr)\ln\frac{\psi_n(t)}{\psi_n(t)+r_n}+r_n\Bigr)\,\dd t,
\end{align*}
where $Y_j^{(n)}$ is the centered version of the process $X_j^{(n)}$.

Since the stochastic integrals with respect to $Y_j^{(n)}$, $j=1,\ldots,n$,
are independent and has mean zero, we obtain
\begin{align*}
\Ex_\vartheta^{(n)}\bigl(\ln Z_{n,\vartheta}(u_1)-\ln
Z_{n,\vartheta}(u_2)\bigr)^2 &=n\,\Ex_\vartheta^{(n)}
\left(\int_{\vartheta+u_1\varphi_n}^{\vartheta+u_2\varphi_n}\ln\frac{\psi_n(t)}{\psi_n(t)+r_n}\,\dd
Y_j^{(n)}(t)\right)^2\\
&\ \quad+n^2\,\biggl(\int_{\vartheta+u_1\varphi_n}^{\vartheta+u_2\varphi_n}
\Bigl(\bigl(\psi_n(t)+r_n\bigr)\ln\frac{\psi_n(t)}{\psi_n(t)+r_n}+r_n\Bigr)\,\dd
t\biggr)^2\\
&=E_1+E_2
\end{align*}
with obvious notations.

Using elementary inequalities $\ln(1+x)\leq x$ and $\ln(1+x)\geq x-x^2/2$ for
$\left|x\right|<1/2$, for sufficiently large values of $n$ (such that
$\frac{r_n}{\psi_n(t)+r_n}<\frac{r_n}{\ell}<1/2$) we obtain
$$
-\frac{r_n}{\psi_n(t)+r_n}-\frac{r_n^2}{2\bigl(\psi_n(t)+r_n\bigr)^2} \leq \ln
\frac{\psi_n(t)}{\psi_n(t)+r_n} \leq -\frac{r_n}{\psi_n(t)+r_n}\,.
$$

For $E_1$, if $r_n\leq 0$, we obtain
\begin{align*}
E_1 &=n\int_{\vartheta+u_1\varphi_n}^{\vartheta+u_2\varphi_n}\Bigl(\ln
\frac{\psi_n(t)}{\psi_n(t)+r_n}\Bigr)^2(\psi_n(t)+r_n)\,\dd t\\
&\leq n\int_{\vartheta+u_1\varphi_n}^{\vartheta+u_2\varphi_n}
\frac{r_n^2}{\psi_n(t)+r_n}\;\dd t\leq
n\,\frac{(u_2-u_1)r_n^2\varphi_n}{\ell}=\frac{\left|u_1-u_2\right|}{\ell}\,.
\end{align*}
As to the case $r_n\geq 0$, as $\frac{r_n}{\psi_n(t)+r_n}<1/2$, we have
\begin{align*}
E_1 &=n\int_{\vartheta+u_1\varphi_n}^{\vartheta+u_2\varphi_n}\Bigl(\ln
\frac{\psi_n(t)}{\psi_n(t)+r_n}\Bigr)^2(\psi_n(t)+r_n)\,\dd t\\
&\leq n\int_{\vartheta+u_1\varphi_n}^{\vartheta+u_2\varphi_n}
\left[\frac{r_n^2}{\psi_n(t)+r_n}+ \frac{r_n^3}{(\psi_n(t)+r_n)^2}+
  \frac{r_n^4}{4(\psi_n(t)+r_n)^3}\right]\;\dd t\\
&\leq n\int_{\vartheta+u_1\varphi_n}^{\vartheta+u_2\varphi_n}
\frac{r_n^2}{\psi_n(t)+r_n}\left[1+\frac{1}{2}+\frac{1}{16}\right]\;\dd t
\leq\frac{25\left|u_1-u_2\right|}{16\ell}\,.
\end{align*}

For $E_2$, we have
$$
-\frac{r_n^2}{2\ell} \leq -\frac{r_n^2}{2(\psi_n(t)+r_n)} \leq
\bigl(\psi_n(t)+r_n\bigr)\ln \frac{\psi_n(t)}{\psi_n(t)+r_n}+r_n \leq 0,
$$
and hence
\begin{align*}
E_2&=n^2\,\biggl(\int_{\vartheta+u_1\varphi_n}^{\vartheta+u_2\varphi_n}\,\Bigl(\bigl(\psi_n(t)+r_n\bigr)\ln
\frac{\psi_n(t)}{\psi_n(t)+r_n}+r_n\Bigr)\,\dd t\biggr)^2\\
&\leq n^2\biggl(\int_{\vartheta+u_1\varphi_n}^{\vartheta+u_2\varphi_n}\,\frac{r_n^2}{2\ell}\;\dd t\biggr)^2
=\frac{(u_2-u_1)^2}{4\ell^2}\,.
\end{align*}

Thus, for sufficiently large values of $n$, we have
$$
\Ex_\vartheta^{(n)}\bigl(\ln Z_{n,\vartheta}(u_1)-\ln Z_{n,\vartheta}(u_2)\bigr)^2\leq
\frac{25\left|u_1-u_2\right|}{16\ell}+\frac{(u_2-u_1)^2}{4\ell^2}\,.
$$

In the case $u_1,u_2\leq 0$, proceeding similarly, we obtain the same
inequality.

Finally, in the case $u_1 u_2<0$ (say $u_1<0$ and $u_2>0$), we obtain
\begin{align*}
\Ex_\vartheta^{(n)}\bigl(\ln Z_{n,\vartheta}(u_1)-\ln Z_{n,\vartheta}(u_2)\bigr)^2
&\leq 2\,\Ex_\vartheta^{(n)}\bigl(\ln Z_{n,\vartheta}(u_1)\bigr)^2+
2\,\Ex_\vartheta^{(n)}\bigl(\ln Z_{n,\vartheta}(u_2)\bigr)^2\\
&\leq \frac{25\left|u_1\right|}{8\ell}+\frac{u_1^2}{2\ell^2}+
 \frac{25\left|u_2\right|}{8\ell}+\frac{u_2^2}{2\ell^2}\\
&=\frac{25}{8\ell}\bigl(\left|u_1\right|+\left|u_2\right|\bigr)
+\frac{1}{2\ell^2}\bigl(u_1^2+u_2^2\bigr)\\
&\leq\frac{25\left|u_1-u_2\right|}{8\ell}
+\frac{(u_2-u_1)^2}{\ell^2}\,.
\end{align*}

Note that this final inequality holds for all the three cases, and
so
$$
\Pb_\vartheta^{(n)}\bigl(\left|\ln Z_{n,\vartheta}(u_1)-\ln Z_{n,\vartheta}(u_2)\right|
>\varepsilon\bigr)\leq\frac{25\left|u_1-u_2\right|}{8\varepsilon^2\ell}
+\frac{(u_2-u_1)^2}{\varepsilon^2\ell^2}
$$
for all $u_1,u_2\in U_n$ and sufficiently large values of $n$. Hence,
$$
 \lim_{n \to +\infty}\ \sup_{\left|u_1-u_2\right|<h}
 \Pb_\vartheta^{(n)}\bigl(\left|\ln Z_{n,\vartheta}(u_1)-\ln
 Z_{n,\vartheta}(u_2)\right| >\varepsilon\bigr)\leq
 \frac{25h}{8\varepsilon^2\ell}+\frac{h^2}{\varepsilon^2\ell^2}\to 0
$$
as $h\to 0$, and so, the lemma is proved.
\end{lemproof}

\begin{lemproof}{to_0_L6}
We have
$$
\Pb_\vartheta^{(n)}\biggl(\sup_{\left|u\right|>D} Z_{n,\vartheta}(u) >e^{-b
  D}\biggr)\leq \Pb_\vartheta^{(n)}\biggl(\sup_{u>D} Z_{n,\vartheta}(u) >e^{-b
  D}\biggr)+ \Pb_\vartheta^{(n)}\biggl(\sup_{u<-D} Z_{n,\vartheta}(u) >e^{-b
  D}\biggr).
$$

In order to estimate the first term, first let us note that the Markov process
$Z_{n,\vartheta}(u)$, $u\geq 0$, is a martingale.  Indeed, for any $v\geq
u\geq 0$, using the representation \eqref{explicitSI1} we can write
\begin{align*}
\Ex\bigl(&Z_{n,\vartheta}(v)\bigm|Z_{n,\vartheta}(u)\bigr)\\
&=\Ex\biggl(\exp
\biggl\{\sum_{j=1}^n\int_{\left(\vartheta+u\varphi_n,\vartheta+v\varphi_n\right]}\ln
\frac{\psi_n(t)}{\psi_n(t)+r_n}\;X^{(n)}_j(\dd t) +\frac{v-u}{r_n}\biggr\}\,
Z_{n,\vartheta}(u)\biggm|Z_{n,\vartheta}(u)\Biggr)\\
&=Z_{n,\vartheta}(u)\:\exp\Bigl\{\frac{v-u}{r_n}\Bigr\}\,
\prod_{j=1}^n\Ex\exp\biggl\{\int_{\left(\vartheta+u\varphi_n,\vartheta+v\varphi_n\right]}\ln
  \frac{\psi_n(t)}{\psi_n(t)+r_n}\;X^{(n)}_j(\dd t)\biggr\}\\
&=Z_{n,\vartheta}(u)\:\exp\Bigl\{\frac{v-u}{r_n}\Bigr\}\,
  \exp\biggl\{n\int_{\vartheta+u\varphi_n}^{\vartheta+v\varphi_n}
  \Bigl(\frac{\psi_n(t)}{\psi_n(t)+r_n}-1\Bigr)({\psi_n(t)+r_n})\;\dd
  t\biggr\}\\
&=Z_{n,\vartheta}(u)\:\exp\Bigl\{\frac{v-u}{r_n}\Bigr\}\,\exp\biggl\{-n
  \int_{\vartheta+u\varphi_n}^{\vartheta+v\varphi_n} r_n\;\dd t\biggr\}\\
&=Z_{n,\vartheta}(u)\:\exp\Bigl\{\frac{v-u}{r_n}\Bigr\}\,\exp\bigl\{-n
  (v-u)\varphi_n r_n\bigr\}=Z_{n,\vartheta}(u).
\end{align*}

Hence, the process $X(t)=Z_{n,\vartheta}^{1/2}(t+D)$, $t\geq 0$, is a
supermartingale and, using the maximal inequality for positive
supermartingales $\bigl($see, for example, Revuz and
Yor~\cite[Exercise~2.1.15]{RY99}$\bigr)$, we get
$$
\Pb_\vartheta^{(n)}\biggl(\sup_{u>D} Z_{n,\vartheta}(u) > e^{-b D}\biggr)
=\Pb_\vartheta^{(n)}\biggl(\sup_{t>0} X(t) > e^{-b D/2}\biggr)\leq e^{b
  D/2}\,\Ex X(0)=e^{b D/2}\,\Ex Z_{n,\vartheta}^{1/2}(D).
$$

So, using Lemma~\ref{L4}, we can majorate the first term as
$$
\Pb_\vartheta^{(n)}\biggl(\sup_{u>D} Z_{n,\vartheta}(u) > e^{-b D}\biggr)\leq
e^{b D/2}\,e^{-k_* D}\leq e^{-b D},
$$
where the last inequality is valid if $b\leq 2k_*/3$.

For the second term, in a similar manner (and under the same condition $b\leq
2k_*/3$) we obtain the bound
$$
\Pb_\vartheta^{(n)}\biggl(\sup_{u<-D} Z_{n,\vartheta}(u) > e^{-b D}\biggr)\leq
e^{-b D},
$$
and so, the required inequality holds with $C=2$ and any $b\in(0,2k_*/3]$.
\end{lemproof}

\begin{acknowledgements}
This study was partially supported by the Russian Science Foundation (research
project no.~14-49-00079).  The authors would like to thank the anonymous
reviewers for their helpful and constructive comments, which greatly
contributed to the improvement of the paper.
\end{acknowledgements}

\end{document}